\documentclass[12pt]{amsart}
\usepackage{amsmath}
\usepackage{amsthm}
\usepackage{amscd}
\usepackage{amsfonts}
\usepackage{amssymb}
\usepackage{latexsym}

\newtheorem{thm}{Theorem}[section]
\newtheorem{lem}[thm]{Lemma}
\newtheorem{theo}{Theorem}

\newtheorem{cor}[thm]{Corollary}
\newtheorem{prop}[thm]{Proposition}

\theoremstyle{definition}

\newtheorem{defn}[thm]{Definition}

\newcommand{\seq}{{\lesssim}} 

\newcommand{\ev}{\mathrm{ev}}
\newcommand{\res}{\mathrm{res}}
\newcommand{\ind}{{\mathrm{ind}}}
\newcommand{\map}{{\mathrm{map}}}
\newcommand{\Ocal}{{\mathcal{O}}}
\newcommand{\uO}{{\underline{\mathcal{O}}}}

\newcommand{\W}{{\mathbb{W}}}
\newcommand{\IW}{{\mathbb{I}}}
\newcommand{\Z}{{\mathbb{Z}}}
\newcommand{\C}{{\mathcal{C}}}
\newcommand{\N}{{\mathbb{N}}}
\newcommand{\op}{{\operatorname{op}}}
\newcommand{\pr}{\operatorname{pr}}
\newcommand{\ob}{\operatorname{ob}}
\newcommand{\id}{\operatorname{id}}

\newcommand{\Set}{\mathcal{E}\! \mathit{ns}}
\newcommand{\Ens}{\mathcal{F}\!\mathit{in}}
\newcommand{\colim}[1]{\mathop{\underset{#1}
            {{\text{\rm colim}}}}}

\author{M. Brun}
\title{ Witt vectors and Tambara Functors}
\begin{document}
\begin{abstract}
We give, for every finite group $G$, a combinatorial
description of the 
ring of $G$-Witt vectors on a polynomial algebra over the
integers. Using this description we show that the functor, which takes a
commutative ring with trivial action of $G$ to its ring of Witt vectors,
coincides 
with the left adjoint of the algebraic functor from the category of
$G$-Tambara functors to the category of commutative rings with an action
of $G$. 
\end{abstract}
\email{brun@mathematik.uni-osnabrueck.de}
\address{FB Mathematik/Informatik, Universit\"at Osnabr\"uck, 49069
Osnabr\"uck, Germany}
\date{April 30, 2003}
\thanks{The author wishes to thank the University of
Osnabr\"uck for hospitality and Rainer Vogt for a number of helpful
conversations. This work has been supported by the DFG}
\maketitle
\section*{Introduction}
In \cite{Witt} Witt constructed an endofunctor on the category of
commutative rings, which takes a commutative ring $R$ to the ring
$\W_p(R)$ of $p$-typical Witt vectors. This construction can be used to
construct field extensions of the $p$-adic numbers \cite{Schmid-Witt},
and it is essential 
in the construction of crystalline 
cohomology \cite{Berthelot}. In \cite{DS} Dress and Siebeneicher
constructed for every 
pro-finite group $G$ an endofunctor $\W_G$ on the category of
commutative rings. In the case where $G = \widehat C_p$ is
the pro-$p$-completion of the infinite cyclic group
the functors $\W_p$ and $\W_G$ agree. The functor
$\W_G$ is constructed in such a way that $\W_G(\Z)$ is an appropriately
completed 
Burnside ring for the pro-finite group $G$. For an arbitrary commutative
ring $R$,
the ring $\W_G (R)$ is somewhat mysterious, even when $G = \widehat C_p$. 
The first aim of the present paper is to give a new description of
the ring $\W_G(R)$ when $R$ is a polynomial algebra over the integers. 
In the special case where $G$ is finite our description is given in terms of a
category $U^G$ introduced by 
Tambara in \cite{Tambara}. The category
$U^G$ has as objects the class 
of finite $G$-sets, and the morphism set $U^G(X,Y)$ is constructed as follows:
Let $U_+^G(X,Y)$ denote the set of equivalence
classes of chains of
$G$-maps of the form $X \leftarrow A \to B \to Y$, where two
chains $X \leftarrow A \to B \to Y$ and 
$X \leftarrow A' \to B' \to Y$ are defined to be equivalent if there
exist $G$-bijections $A \to A'$ and $B \to B'$ making
the following diagram commutative:
\begin{displaymath}
  \begin{CD}
    X @<<< A @>>> B @>>> Y \\
    @| @VVV @VVV @| \\
    X @<<< A' @>>> B' @>>> Y.
  \end{CD}
\end{displaymath}
The set $U_+^G(X,Y)$ has a natural semi-ring-structure described in
\ref{UGprod}. By group-completing the underlying additive monoid of
$U^G_+(X,Y)$ we obtain a commutative ring, and we define $U^G(X,Y)$ to
be the set obtained by forgetting the ring-structure.
For a subgroup $H$ of $G$ we let $G/H$ denote the $G$-set of left cosets
of$H$ in $G$, and we let $G/e$ denote $G$ considered as a left $G$-set. 
The following theorem is a special case of theorem \ref{subwitt}. 
\begin{theo}
  Let $G$ be a finite group and let $X$ be a finite set considered as a
  $G$-set with trivial action. The 
  ring $U^G(X,G/e)$ is the polynomial ring 
  $\Z[X]$ over $\Z$, with one indeterminate for each element in $X$, and the
  ring $U^G(X,G/G)$ is naturally isomorphic to the ring
  $\W_G(U^G(X,G/e)) = \W_G(\Z[X])$.
\end{theo}
Our second aim is to advertise the category of $G$-Tambara functors,
that is, the category $[U^G,\Set]_0$ of
set-valued functors on $U^G$ preserving finite products. Tambara calls
such a functor a $\mathrm{TNR}$-functor, an acronym for ``functor with
trace, norm and restriction'', but we have chosen to call them  
$G$-Tambara functors. This category is intimately related to the Witt
vectors of Dress and Siebeneicher. In order to explain this relation
we consider the full subcategory $U^{fG}$ of $U^G$ with free $G$-sets
as objects, and we note that the category of $fG$-Tambara functors, that
is, the category $[U^{fG},\Set]_0$ of
set-valued functors on $U^{fG}$ preserving finite products is
equivalent to the category of commutative rings with an action of $G$
through ring-automorphisms. 
The categories $U^G$ and $U^{fG}$ are coloured theories in the sense of
Boardman and Vogt \cite{BV}, and therefore they are complete and
cocomplete, and
the forgetful functor $[U^{fG},\Set]_0 \to
[U^{G},\Set]_0$ induced by the inclusion $U^{fG} \subseteq U^G$ has a
left adjoint functor, which we shall denote by $L_G$.  
\begin{theo}\label{mainwitt0}
  Let $G$ be a finite group. If $R$ is a commutative ring with $G$ acting
  trivially, then for every subgroup 
  $H$ of $G$, there is a natural
  isomorphism $\W_H(R) \cong L_G R(G/H)$.
\end{theo}
We are also able to describe the ring $L_GR(G/H)$ in the case where $G$
acts nontrivially on $R$. In order to do this we need to introduce an
ideal ${\IW}_G(R)$ of $\W_G(R) = R^{\uO(G)}$, where $\uO(G)$ is the set of
conjygacy classes of subgroups of $G$. An element $a \in \W_G(R)$ is of
the form $a = (a_H)'_{H \leq G}$, where the prime means that $a_H =
a_{gHg^{-1}} \in R$ for $H \leq G$ and $g \in G$.
We let ${\IW}_G(R)
\subseteq \W_G(R)$ denote the ideal 
generated by elements of the form $a-b$, where $a=(a_K)'_{K \leq
G}$ and $b = (b_K)'_{K \leq G}$ satisfy that for every $K \leq G$
there exist $n \ge 1$ and $g_{1,K}, \dots g_{n,K} \in N_K(G)$ and
$a_{1,K},\dots,a_{n,K} \in R$ such that firstly $g_{1,K}K=  \dots =
g_{n,K}K$ 
and secondly
$a_K = a_{1,K} \dots a_{n,K}$ and $b_K = (g_{1,K} a_{1,K}) \dots
(g_{n,K} a_{n,K}$). 
Our description of $L_GR(G/H)$ is formulated in the following
generalization of theorem \ref{mainwitt0}:
\begin{theo}\label{mainwitt1}
  Let $G$ be a finite group and let $R$ be a commutative ring with an
  action of $G$.
  There is a natural isomorphism $ 
  \W_H(R)/{\IW}_H(R) \cong
  (L_GR)(G/H)$. 
\end{theo}

Let us stress that in our setup the above result only makes sense for a
finite group. For a
pro-finite group $G$ it might be possible to construct 
appropriate variants $\widehat U^G(X,Y)$ of the rings $U^G(X,Y)$ such
that the  the ring $\widehat U^G(X,G/G)$ is naturally isomorphic to
$\W_G(\Z[X])$. We shall not go further into this here.

There have been given other combinatorial
constructions of the Witt vectors by Metropolis and Rota \cite{MR} and
Graham \cite{Graham}. The ring of Witt vectors has been related to the
ring of necklaces by Dress and Siebeneicher in \cite{DS2}.
Our combinatorial description differs from the previous ones mostly in
that it is forced by some additional structure on the Witt vectors.

It is our belief that the category of $G$-Tambara functors has
applications yet to be discovered. There is a rich supply
of them coming from equivariant stable homotopy theory. In fact, every
$E_\infty$ ring $G$-spectrum gives rise to a $G$-Tambara functor by
taking the zeroth homotopy group \cite{tamfun2}. 
In the case where $G = e$ is the trivial group, the category of
$G$-Tambara functors is equivalent to the category of commutative
rings. It is well known that every commutative ring can be realized as
the zeroth homotopy group of an $E_\infty$ ring-spectrum.
For an arbitrary group $G$ one may speculate whether  
$G$-Tambara functors can be realized as the zeroth homotopy group of an
$E_\infty$ ring-spectrum with an action of $G$. 

Some concepts from commutative algebra can be carried over to the setting of
$G$-Tambara functors. Most notably, the notions of an ideals,
modules and chain complexes have $G$-Tambara versions. 
We leave these concepts to future work.

The paper is organized as follows:
In section \ref{prereq} we have collected some of the results from the papers
\cite{DS} and \cite{Tambara} that we need in the rest of the paper. In
section \ref{colth} we note that the category of Tambara functors is the
category of algebras for a coloured theory. In
section \ref{tmap} we construct a homomorphism relating Witt vectors
and Tambara 
functors, which we have chosen to call the Teichm\"uller
homomorphism. In section \ref{wipoly} we prove the fundamental fact that the
Teichm\"uller homomorphism is a ring-homomorphism. In section
\ref{freeta} we prove 
that for free Tambara functors, the Teichm\"uller homomorphism is an
isomorphism, and finally in section \ref{witafun} we prove theorem
\ref{mainwitt1}.

\section{Prerequisites}
\label{prereq}
In this section we fix some notation and recollect results from
\cite{DS} and  \cite{Tambara}.
All rings are supposed to be both commutative and unital. Given a group
$G$ we only consider left actions of $G$. A $G$-ring is
a ring with an action of the group $G$ through ring-automorphisms.

Given a pro-finite group $G$ we let $\Ocal(G)$ denote the $G$-set of
open  
subgroups of $G$ with action given by conjugation and
we let
${\uO}(G)$ denote the set of 
conjugacy classes of open subgroups of $G$. For a $G$-set $X$ and a subgroup
$H$ of $G$ we define 
$|X^H|$ to be the cardinality of the set $X^H$ of $H$-invariant
elements of $X$.
The following is
the main result
of \cite{DS}.
\begin{thm}\label{DSmain}
  Let $G$ be a pro-finite group.
  There exists a unique endofunctor $\W_G$ on the category of
  rings such that for a ring $R$ the ring $\W_G(R)$
  has 
  the set $R^{\uO(G)}$ of maps from the set
  ${\uO(G)}$ to $R$ as underlying set, in such a way that for
  every ring-homomorphism $h : R \to R'$ and every $x \in \W_G(R)$ one has
  $\W_G(h)(x) = h \circ x$, while for any subgroup $U$ of $G$ the family
  of $G$-maps
  \begin{displaymath}
    \phi^R_U : \W_G(R) \to R
  \end{displaymath}
  defined by 
  \begin{displaymath}
    x = (x_V)'_{V \leq G}  \mapsto {\sum}'_{U \seq V \leq G} |U(G/V)^U| \cdot
    x_V^{(V:U)} 
  \end{displaymath}
  provides a natural transformation from the functor $\W_G$ into the
  identity functor. Here $U \seq V$ means that the subgroup $U$ of $G$ is
  subconjugate to $V$, i.e., there exists some $g \in G$ with $U
  \leq gVg^{-1}$, $(V:U)$ means the index of $U$ in $gVG^{-1}$
  which coincides with $(G:U)/(G:V)$ and therefore is independent of
  $g$, and the symbol ``$\sum'$'' is meant to indicate that for each
  conjugacy class of subgroups $V$ with $U \seq V$ exactly one summand
  is taken, and 
  an element $a \in \W_G(R)$ is written on the form $a = (a_V)'_{V
  \leq G}$, where the prime means that $a_V = a_{gVg^{-1}}$ for $g \in G$.
\end{thm}
In the section \ref{tmap} we shall reprove
theorem \ref{DSmain}.

The rest of this section is a recollection of the work \cite{Tambara} of
Tambara. 
We let $\Ens$ denote the category
of finite sets and we let $\Ens^G$ denote the category of finite
$G$-sets.
Given a finite $G$-set $X$ we denote by $\Ens^G/X$
the category of objects over $X$ in $\Ens^G$.
Given $f : X \to Y$ in $\Ens^G$ the pull-back functor
\begin{eqnarray*}
  \Ens^G/Y &\to& \Ens^G/X \\
  (B \to Y) &\mapsto& (X \times_Y B \to X)
\end{eqnarray*}
has a right adjoint
\begin{eqnarray*}
  \Pi_f : \Ens^U/X &\to& \Ens^G/Y \\
  (A \xrightarrow p X) &\mapsto& (\Pi_f A \xrightarrow{\Pi_f p} Y),
\end{eqnarray*}
where $\Pi_f p$ is made from $p$ as follows. For each $y \in Y$, the
fibre $(\Pi_f p)^{-1}(y)$ is the set of maps $ s : f^{-1}(y) \to A$ such
that $p(s(x)) = x$ for all $x \in f^{-1}(y)$. If $g \in G$ and $s \in
(\Pi_fp)^{-1}(y)$, the map $^g s : f^{-1}(g y) \to A$ taking $x$ to $g
s(g^{-1} x)$ belongs to $(\Pi_f p)^{-1}(gy)$. The operation $(g,s)
\mapsto ^gs$ makes $\Pi_f A$ a $G$-set and $\Pi_f p$ a $G$-map.

There is a commutative diagram
\begin{displaymath}
\begin{CD}
  X @<p<< A @<e<< X \times_Y \Pi_f A \\
  @VfVV @. @Vf'VV \\
  Y @<{\Pi_f p}<< \Pi_f A @= \Pi_f A,
\end{CD}
\end{displaymath}
where $f'$ is the projection and $e$ is the evaluation map $(x,s)
\mapsto s(x)$. A diagram in $\Ens^G$ which is isomorphic to a diagram of
the above form is called an {\em exponential diagram}.

We say that two diagrams $X \leftarrow A \to B \to Y$ and $X \leftarrow
A' \to B' \to Y$ in $\Ens^G$ are equivalent if there exist
$G$-isomorphisms $A \to A'$, $B \to B'$ making the diagram
\begin{displaymath}
\begin{CD}
  X @<<< A @>>> B @>>> Y \\
  @| @VVV @VVV @| \\
  X @<<< A' @>>> B' @>>> Y   
\end{CD}
\end{displaymath}
commutative, and we let $U^G_+(X,Y)$ be the set of the equivalence classes
$[X \leftarrow A \to B \to Y]$ of diagrams $X \leftarrow A \to B \to
Y$. 

Tambara defines an operation $\circ : U^G_+(Y,Z) \times U^G_+(X,Y)
\to U^G_+(X,Z)$ by
\begin{displaymath}
  [Y \leftarrow C \to D \to Z] \circ [X \leftarrow A \to B \to Y] = [X
  \leftarrow A'' \to \widetilde D \to Z],  
\end{displaymath}
where the maps on the right are composites of the maps in the
diagram 
\begin{displaymath}
\begin{CD}
  X @<<< A @<<< A' @<<< A'' \\
  @. @VVV @VVV @VVV \\
  @. B @<<< B' @<<< \widetilde C \\
  @. @VVV @VVV @| \\
  @. Y @<<< C @. \widetilde C \\
  @. @. @VVV @VVV \\
  @. Z @<<< D @<<< \widetilde D.
\end{CD}
\end{displaymath}
Here the three squares are pull-back diagrams and the diagram
\begin{displaymath}
\begin{CD}
  C @<<< B' @<<< \widetilde C \\
  @VVV @. @VVV \\
  D @<<< \widetilde D @= \widetilde D
\end{CD}
\end{displaymath}
is an exponential diagram.
He verifies that $U^G_+$ is a category with $\circ$ as commposition
and given $f : X \to Y$ in $\Ens^G$ he introduces the notation
\begin{eqnarray*}
  R_f &=& [Y \xleftarrow f X \xrightarrow = X \xrightarrow = X], \\
  T_f &=& [X \xleftarrow = X \xrightarrow = X \xrightarrow f Y]
  \qquad \text{and}\\
  N_f &=& [X \xleftarrow = X \xrightarrow f Y \xrightarrow = Y].
\end{eqnarray*}
Every morphism in $U^G_+$ is a composition of
morphisms on the above form.
He also shows that every
object $X$ of $U^G_+$ is a semi-ring object in the following natural way:
\begin{prop}\label{UGprod}
  Given objects $X$ and $Y$ in $U^G_+$, there is semi-ring-structure on
  $U_+^G(X,Y)$ given as follows: 
  \begin{eqnarray*}
  0 &=& [X \leftarrow \emptyset \to \emptyset \to Y], \\
  1 &=& [X \leftarrow \emptyset \to Y \to Y], \\
  {[X \leftarrow A \to B \to Y] + [X \leftarrow A' \to B' \to
  Y]} 
  &=& [X \leftarrow A \amalg A' \to B \amalg B' \to Y] 
  \end{eqnarray*}
  and
  \begin{multline*}
  {[X \leftarrow A \to B \to Y] \cdot [X \leftarrow A' \to B' \to
  Y]} 
  = \\ [X \leftarrow B \times_Y A' \amalg A \times_Y B' \to B \times_Y 
  B' \to Y].
  \end{multline*}
\end{prop}
It is also shown in \cite{Tambara} that there is a unique category $U^G$
satisfying the following conditions:
\begin{enumerate}
  \item[(i)] $\ob U^G = \ob U^G_+$ 
  \item[(ii)] The morphism set $U^G(X,Y)$ is the group completion of the
  underlying additive monoid of $U^G_+(X,Y)$.
  \item[(iii)] The group completion maps $k : U^G_+(X,Y) \to U^G(X,Y)$ and the
  identity on $\ob(U^G_+)$ form a functor $k : U^G_+ \to U^G$.
  \item[(iv)] The functor $k$ preserves finite products.
\end{enumerate}
\begin{prop}\label{semi-ring-str}
(i) If $X_1 \xrightarrow {i_1} X \xleftarrow {i_2} X_2$ is a sum diagram
in $\Ens^G$, then 
$X_1 \xleftarrow {R_{i_1}} X \xrightarrow {R_{i_2}} X_2$ is a product
diagram in $U^G$ and $\emptyset$ is final in $U^G$.

(ii) Let $X$ be a $G$-set and $\nabla : X \amalg X \to X$ the folding
map, $i : \emptyset \to X$ the unique map. Then $X$ has a structure of a
ring object of $U^G$ with addition $T_{\nabla}$, additive unit $T_i$,
multiplication $N_{\nabla}$ and multiplicative unit $N_i$.

(iii) If $f : X \to Y$ is a $G$-map, then the morphisms $R_f$, $T_f$ and
$N_f$ of $U^G$ preserve the above structures of ring, additive group and
multiplicative monoid on $X$ and $Y$, respectively.
\end{prop}
Given a category $\C$ with finite products, we shall denote the category
of set-valued functors on $\C$ preserving finite products by
$[\C,\Set]_0$. The morphisms in $[\C,\Set]_0$ are given by natural
transformations. 
\begin{defn}\label{tamdef}
  The category of $G$-Tambara functors is the category
  $[U^G,\Set]_0$. That is, 
  a $G$-Tambara functor (called $\mathrm{TNR}$-functor in
  \cite{Tambara}) is a set-valued functor on
  $U^G$ preserving finite products, and a morphism of $G$-Tambara functors
  is a natural 
  transformation. 
\end{defn}
Given a $G$-Tambara functor $S$ and $[X \leftarrow A \to B \to Y] \in
U^G(X,Y)$ we obtain a function $S[X \leftarrow A
  \to B \to Y] : S(X) \to S(Y)$. Since $S$ is product-preserving,
it follows from (ii) of proposition \ref{semi-ring-str} that $S(X)$ is a
ring. 
Given a finite $G$-map $f : X \to Y$
we shall use the notation $S^*(f) = S(R_f)$, $S_+(f) = S(T_f)$ and
$S_\bullet(f) = S(N_f)$.
It follows from (iii) of \ref{semi-ring-str} that $S^*(f)$ is a
ring-homomorphism, that 
$S_+(f)$ is an additive homomorphism and that $S_\bullet(f)$ is
multiplicative. A $G$-Tamara functor $S$ is uniquely determined by the
functions $S^*(f)$, $S_+(f)$ and $S_\bullet(f)$ for all $f : X \to Y$ in
$\Ens^G$. 

Given subgroups $K \leq H \leq G$ we shall denote by $\pi^H_K : G/K
\to G/H$ the projection induced by the inclusion $K \leq H$ and given
$g \in G$ we shall let $c_g : G/H \to G/gHg^{-1}$ denote conjugation
by $g$, $c_g(\sigma H) = \sigma g^{-1} (g H g^{-1})$.

\section{Coloured theories}
\label{colth}
In this section we shall explain that the category $U^G$ is an
$\uO(G)$-coloured 
category in the sense of Boardman and Vogt \cite{BV}.
\begin{defn}[{\cite[2.3,]{BV}}]
  \begin{enumerate}
  \item[(i)]
  Let $\Ocal$ be a finite set. An {\em $\Ocal$-coloured theory} is a category
  $\Theta$ together with a faithful functor $\sigma_{\Theta} : (\Ens
  /\Ocal)^{\op} \to \Theta$ such that firstly $\sigma_{\Theta}$ preserves
  finite products and secondly every object of $\Theta$ is isomorphic to
  an object in the image of $\sigma_{\Theta}$.
  \item[(ii)]
  The {\em category of algebras} over a theory $\Theta$ is the category
  $[\Theta, \Set]_0$ of product-preserving set-valued functors on
  $\Theta$. 
  \end{enumerate}
  \item[(iii)]
  A morphism $\gamma : \Theta \to \Psi$ of $\Ocal$-coloured theories is a
  functor preserving finite products.
\end{defn}
Other authors, e.g. \cite{AR} and \cite{Rezk}, use the name ``sorted
theory'' for a coloured theory. 

Given a finite group $G$, choosing representatives $G/H$ for the objects
of $\uO(G)$, we can construct a functor
\begin{eqnarray*}
  \sigma_{({\Ens^G})^{\op}} : (\Ens /\uO(G))^{\op} &\to& (\Ens^G)^{\op}, \\
  (z : Z \to \uO(G)) &\mapsto& \coprod_{[G/H] \in \uO(G)} G/H
  \times z^{-1}([G/H]).
\end{eqnarray*}
This way we give $(\Ens^G)^{\op}$ the structure of an $\uO(G)$-coloured
theory. Composing $\sigma_{(\Ens^G)^{\op}}$ with the functor $R :
(\Ens^G)^{\op} \to U^G$, $f \mapsto R_f$ we obtain by (i) of
\ref{semi-ring-str} a 
functor $\sigma_{U^G} : (\Ens/\uO(G))^{\op} \to U^G$ preserving finite
products, and we obtain a structure of $\uO(G)$-coloured theory on
$U^G$. The functor $R$ is a morphism of $\uO(G)$-coloured theories.

Let $V^G \subseteq U^G$ denote the subcategory of $U^G$ with the same
class of objects as $U^G$ and with $V^G(X,Y) \subseteq U^G(X,Y)$ given
by the subgroup generated by morphisms of the form $[X \leftarrow A
\xrightarrow = A \to Y]$. The inclusion $V^G \subseteq U^G$ preserves
finite products, and we have morphisms $(\Ens^G)^{\op} \to V^G \to U^G$
of $\uO(G)$-coloured theories. The category $V^G$ is strongly related
to the category spans considered by Lindner in \cite{Lindner}, and in
fact the
category of Mackey functors is equal to the category of
$V^G$-algebras. Let $\Ens^{fG}$ denote the full subcategory of $\Ens^G$
with finite free $G$-sets as objects. The functor
\begin{displaymath}
  \sigma_{(\Ens^{fG})^{\op}} : (\Ens /\uO(G))^{\op} \to (\Ens^{fG})^{\op},
  \qquad (z : Z \to \uO(G)) \mapsto G/e
  \times z^{-1}([G/e])
\end{displaymath}
gives $(\Ens^{fG})^{\op}$ the structure of an $\uO(G)$-coloured
theory. Similarly we can consider the $\uO(G)$-coloured theories given
by the full subcategories $U^{fG} \subseteq U^G$ and $V^{fG} \subseteq
V^G$ with finite free $G$-sets as objects. We have the following diagram
of morphisms of $\uO(G)$-coloured theories:
\begin{displaymath}
\begin{CD}
  (\Ens^{fG})^{\op} @>>> V^{fG} @>>> U^{fG} \\
  @VVV @VVV @VVV \\
  (\Ens^{G})^{\op} @>>> V^{G} @>>> U^{G},
\end{CD}
\end{displaymath}
where the vertical functors are inclusions of full subcategories. 
\begin{lem}\label{Gringistam}
  \begin{enumerate}
  \item[(i)] The category $[U^{fG},\Set]_0$ of $fG$-Tambara functors is
  equivalent to the category of $G$-rings.
  \item[(ii)] The category $[V^{fG},\Set]_0$ of $V^{fG}$-algebras is
  equivalent to the category of left $\Z[G]$-modules.
  \item[(iii)] The category $[(\Ens^{fG})^{\op},\Set]_0$ of
  $(\Ens^{fG})^{\op}$-algebras is 
  equivalent to the category of $G$-sets.
  \end{enumerate}
\end{lem}
\begin{proof}
  Since the statements have similar proofs we only give the proof of (i). 
  Given an $fG$-Tambara functor $S$, we construct a $G$-ring-structure
  on $R = S(G/e)$. Indeed by (ii) of \ref{semi-ring-str} $S(G/e)$ is a
  ring, and 
  given $g \in G$ the right multiplication $g : G/e \to G/e$, $x \mapsto
  xg$, induces a ring-automorphism $S^*(g^{-1})$ of $R = S(G/e)$. From
  the functoriality of $S$ we obtain that $R$ is a $G$-ring. Conversely,
  given a $G$-ring $R'$ we shall construct an $fG$-Tambara functor
  $S'$. We define $S'(X)$ to be the set of $G$-maps from $X$ to $R'$. Given
  $[X \xleftarrow d A \xrightarrow f B \xrightarrow g Y] \in U^G(X,Y)$,
  we define $S'[X \xleftarrow d A \xrightarrow f B \xrightarrow g Y] :
  S'(X) \to S'(Y)$ by the formula
  \begin{displaymath}
    S'[X \xleftarrow d A \xrightarrow f B \xrightarrow g Y](\phi)(y) =
  \sum_{b \in g^{-1}(y)} ( \prod_{a \in f^{-1}(b)} \phi(d(a)))
  \end{displaymath}
  for $\phi \in S'(X)$ and $y \in Y$. We leave it to the reader to check
  that $S \mapsto R$ and $R' \mapsto S'$ are inverse functors up to
  isomorphism. 
\end{proof}
We refer to \cite[Propositions 4.3 and 4.7.]{Rezk} for a proof of the following
two results. Alternatively the reader may modify the proofs given in
\cite[3.4.5 and 3.7.7]{Borceux} for their monocrome versions.
\begin{prop}
  Let $\Theta$ be an $\Ocal$-coloured theory. The category of
  $\Theta$-algebras is complete and cocomplete. 
\end{prop}
\begin{prop}\label{leftKan}
  Given a morphism $\gamma : \Theta \to \Psi$ of $\Ocal$-coloured theories,
  the functor $\gamma^* : [\Psi,\Set]_0 \to [\Theta,\Set]_0$, $A \mapsto
  A \circ \gamma$ has a left adjoint $\gamma_* : [\Theta,\Set]_0 \to
  [\Psi,\Set]_0$. 
\end{prop}
\begin{defn}
  The category of {\em $fG$-Tambara functors} is the category
  $[U^{fG},\Set]_0$ of $U^{fG}$-algebras.
\end{defn}
\begin{defn}
  We let $L_G = i_* : [U^{fG},\Set]_0 \to [U^G,\Set]_0$ denote the left
  adjoint of the functor $i^* : [U^{G},\Set]_0 \to [U^{fG},\Set]_0$
  induced by the inculsion $i : U^{fG} \subseteq U^G$.
\end{defn}
Let us note that $L_G$ can be constructed as the left Kan extension
along $i$, and that for $R \in [U^{fG},\Set]_0$, we have for every
finite free $G$-set $X$ an isomorphism $L_GR(X) \cong R(X)$ because
$U^{fG}$ is a full subcategory of $U^G$.
\section{The Teichm\"uller homomorphism}
\label{tmap}
We shall now
give a connection between the category of $G$-Tambara functors and the
category of rings with an action of $G$.
Recall that $U^{fG}$ denotes the full subcategory of $U^G$ with objects
given by 
free $G$-sets. Below we use notation introduced in \ref{DSmain}. 
\begin{defn}
  Let $G$ be a finite group and let $S$ be a $G$-Tambara
  functor. The {\em unrestricted Teichm\"uller 
  homomorphism} 
  $t : \W_G(S(G/e)) \to S(G/G)$ takes $(x_U)'_{U \leq G} \in
  \W_G(S(G/e))$ to $$t ((x_U)'_{U \leq G}) = {\sum}'_{U\leq G}
  S_+(\pi^G_U) S_\bullet(\pi^U_e)(x_U).$$
\end{defn}
We shall prove the following proposition in the next section.
\begin{prop}\label{tauring}
  The unrestricted Teichm\"uller homomorphism $t : \W_G(S(G/e)) \to
  S(G/G)$ is a ring-homomorphism.  
\end{prop}
In general $t$ will neither be injective nor surjective. 
However, in the introduction, we have given an ideal ${\IW}_G(S(G/e))
\subseteq \W_G(S(G/e))$ contained
in the kernel of $t$:
\begin{defn}\label{PGideal0}
Let $G$ be a finite group and let $R$ be a commutative ring with an
action of $G$. We let ${\IW}_G(R)
\subseteq \W_G(R)$ denote the ideal 
generated by elements of the form $a-b$, where $a=(a_K)'_{K \leq
G}$ and $b = (b_K)'_{K \leq G}$ satisfy that for every $K \leq G$
there exist $n \ge 1$ and $g_{1,K}, \dots g_{n,K} \in N_K(G)$ and
$a_{1,K},\dots,a_{n,K} \in R$ such that firstly $g_{1,K}K=  \dots =
g_{n,K}K$ 
and secondly
$a_K = a_{1,K} \dots a_{n,K}$ and $b_K = (g_{1,K} a_{1,K}) \dots
(g_{n,K} a_{n,K}$). 
\end{defn}
\begin{lem}
  Let $S$ be a $G$-Tambara functor. The unrestricted Teichm\"uller
  homomorphism 
  $t : \W_G(S(G/e)) \to S(G/G)$ maps the ideal ${\IW}_G(S(G/e))$ to zero.
\end{lem}
\begin{proof}
  It suffices to consider the case where $a_U = 0$ for $U \ne K$. 
  We have that 
  \begin{eqnarray*}
    \lefteqn{S_+(\pi^G_K)S_\bullet(\pi^K_e)(g_1 a_1
  \cdots g_n a_n)}\\
  &=& S[\coprod_1^n G/e \xleftarrow{\coprod_{i=1}^n g_i} \coprod_1^n G/e
  \xrightarrow {\pi^K_e \circ \nabla} G/K \xrightarrow {\pi^G_K} G/G]
  (a_1,\dots,a_n) \\ 
  &=& S[\coprod_1^n G/e \xleftarrow{=} \coprod_1^n G/e
  \xrightarrow {\pi^K_e \circ \nabla} G/K \xrightarrow {\pi^G_K} G/G]
  (a_1,\dots,a_n) \\ 
  &=& S_+(\pi^G_K)S_\bullet(\pi^K_e)(a_1
  \cdots a_n),
  \end{eqnarray*}
  where $\nabla : \coprod_1^n G/e \to G/e$ is the fold map.
\end{proof}
\begin{defn}
  Let $G$ be a finite group and let $S$ be a $G$-Tambara
  functor. The {\em Teichm\"uller 
  homomorphism} $\tau : \W_G(S(G/e))/{\IW}_G(S(G/e)) \to S(G/G)$ is the
  ring-homomorphism induced by the unreduced Teichm\"uller homomorphism
  $t$. 
\end{defn}
The following theorem, proved in section \ref{freeta}, implies the case
$H=G$ of theorem \ref{mainwitt1}. 
\begin{thm}\label{mainwitt}
  Let $G$ be a finite group and let $R$ be an $fG$-Tambara
  functor. The Teichm\"uller homomorphism $$\tau 
  : \W_G(L_GR(G/e))/{\IW}_G(L_GR(G/e)) \to
  L_GR(G/G)$$ is an 
  isomorphism. In particular, if $G$ acts trivially on $R$, the Teichm\"uller
  homomorphism is an isomorphism on the form $\tau
  : \W_G(L_GR(G/e)) \to L_GR(G/G)$
\end{thm}
In the above situation there is an isomorphism $L_GR(G/e)
\cong R(G/e)$.
\section{Witt polynomials}
\label{wipoly}
We
shall prove
the following version of a theorem of Dress and Siebeneicher
\cite[Theorem 3.2]{DS}. 
As we shall see below, theorem
\ref{DSmain} and proposition  \ref{tauring} are immediate consequences of
it. 
\begin{thm}\label{maintam} 
  Given a finite group $G$, there exist unique families $(s_U)'_{U \leq
  G}$, $(p_U)'_{U \leq G}$ of integral
  polynomials 
  \begin{displaymath}
    s_U = s_U^G, \quad p_U = p_U^G \in \Z[a_V,b_V \, | \, U \seq V
  \leq G]
  \end{displaymath}
  in two times as many variables $x_V,y_V$ $(U \seq V \leq G)$ as
  there are conjugacy classes of subgroups $V \leq G$ which contain a
  conjugate of $U$ such that for every $G$-Tambara functor $S$ and
  for every $x = {(x_U)}'_{U \leq G}$ and $y= (y_U)'_{U \leq G} \in
  \W_G(S(G/e))$ we have that
  \begin{eqnarray*}
    \tau(x) + \tau(y) &=& \tau((s_U(x_V,y_V \, | \, U \seq V \leq
    G))'_{U \leq G}) \qquad \text{and} \\
    \tau(x) \cdot \tau(y) &=& \tau((p_U(x_V,y_V \, | \, U \seq V \leq
    G))'_{U \leq G}).
  \end{eqnarray*}
  Similarly, there exist
  polynomials $m_U = 
  m_U^G 
  \in \Z[a_V \, | \, U \seq V \leq G]$ such that 
  for every $G$-Tambara functor $S$ and
  for every $x= (x_U)'_{U \leq G} \in
  \W_G(S(G/e))$ we have that
  $-\tau(x) = \tau((m_U(x_V \,
  | \, U \seq V \leq G')_{U \leq G})$. Further, for every subgroup $H$
  of $G$ we have that
  \begin{eqnarray*}
    \phi_H(x) + \phi_H(y) &=& \phi_H((s_U(x_V,y_V \, | \, U \seq V \leq
    G))'_{U \leq G}) \qquad \text{and} \\
    \phi_H(x) \cdot \phi_H(y) &=& \phi_H((p_U(x_V,y_V \, | \, U \seq V \leq
    G))'_{U \leq G}).
  \end{eqnarray*}
  We shall call the polynomials $s_U,p_U$ and $m_U$ the
  {\em Witt polynomials}.
\end{thm}
\begin{proof}[Proof of theorem \ref{DSmain}]
  We first consider the case where $G$ is finite. Given a ring $R$, we
  define operations $+$ and $\cdot$ on $\W_G(R)$ by 
  defining $(a_U)'_{U \leq G} + (b_U)'_{U \leq G} := (s_U(a_V,b_V \, |
  \, U \seq V \le G))'_{U \leq G}$ and $(a_U)'_{U \leq G} \cdot
  (b_U)'_{U \leq G} := (p_U(a_V,b_V \, | 
  \, U \seq V \le G))'_{U \leq G}$. In the case where $R$ has no
  torsion, the map $\phi : \W_G(R) \to {\prod}'_{U \leq G} R$ with
  $U$'th component $\phi_U$ is injective, and hence $\W_G(R)$ is a
  subring of ${\prod}'_{U \leq G} R$. In the case where $R$ has torsion,
  we can choose a surjective ring-homomorphism $R' \to R$ from a torsion
  free ring 
  $R'$. We obtain a surjection $\W_G(R') \to \W_G(R)$ respecting the
  operations $\cdot$ and $+$. Since $\W_G(R')$ is a ring we can
  conclude that $\W_G(R)$ is a ring. Given a surjective homomorphism $\gamma
  : G \to G'$ of finite groups we obtain a ring-homomorphism
  $\mathrm{restr}^G_{G'} : \W_G(R) \to
  \W_{G'}(R)$ with $\mathrm{restr}^G_{G'}((a_U)'_{U \leq G}) =
  ((b_V)')_{V \leq G'})$, where $b_V = a_{\gamma^{-1}(V)}$. (See
  \cite[(3.3.11)]{DS}.) The easiest way to see that
  $\mathrm{restr}^G_{G'}$ is a ring-homomorphism is to note that
  $(\gamma^{-1}H : \gamma^{-1}U) = (H : U)$ and that
  $\phi_{\gamma^{-1}H}(G/\gamma^{-1}U) = 
  \phi_H(G'/H)$. 
  For the case where $G$ is a
  pro-finite group we note that for the ring
  $\W_G(R)$ can be defined to be the limit $\lim_N \W_{G/N}(R)$ taken
  over all finite factor groups $G/N$, with respect maps on 
  the form $\mathrm{restr}^G_{G'}$. 
\end{proof}
\begin{proof}[Proof of proposition \ref{tauring}]
  Proposition \ref{tauring} follows from the first part
  of theorem \ref{maintam} because we use the Witt polynomials to define
  the ring-stucture on the Witt vectors. 
\end{proof}
  For the rest of this section we fix a finite group $G$.
  We shall follow \cite{DS} closely in the proof of theorem
  \ref{maintam}.
  For the uniqueness of the Witt polynomials we consider the
  representable $G$-Tambara functor $\Omega := U^G(\emptyset,-)$. In
  \cite[Theorem 2.12.7]{DS} it is shown that $\tau : {\prod}'_{U \leq G}
  \Z = \W_G(\Omega(G/e)) \to 
  \Omega_G(\Z)$ is a bijection. Hence the Witt polynomials are unique.
  We shall use the following six lemmas to etablish the existence of
  Witt polynomials with the properties required in theorem \ref{maintam}.
\begin{lem}\label{sumlemma}
  For a subset $A$ of $G$, let $U_A := \{ g \in G \, | \, Ag = A \}$
  denote its stabilizer group and let $i_A := |A/U_A|$ denote the
  number of $U_A$-orbits in $A$. If the set $\mathcal{U}(G)$ of
  subsets of $G$ is considered as a $G$-set via $G \times \mathcal{U}(G)
  \to \mathcal{U}(G)$: $(g,A) \mapsto Ag^{-1}$, then for any $s,t \in
  S(G/e)$ one has
  \begin{displaymath}
    S_{\bullet}(\pi^G_e)(s+t) = \sum_{G\cdot A \in G \backslash \mathcal
  U(G)} S_+(\pi^G_{U_A}) S_\bullet(\pi^{U_A}_e)(s^{i_A} \cdot t^{i_{G-A}}).
  \end{displaymath}
\end{lem}
\begin{proof}
  We
  let $i_1,i_2 : G/e \to G/e \amalg G/e$ denote the two natural
  inclusions. We have an exponential diagram 
  \begin{displaymath}
  \begin{CD}
  G/e @<{\nabla}<< G/e \amalg G/e @<d<< G/e \times \mathcal U (G) \\
  @V{\pi^G_e}VV @. @V{\pr}VV \\
  G/G @<<< \mathcal U(G) @= \mathcal U(G),
  \end{CD}
  \end{displaymath}
  where $d(g,A) = i_1(g)$ if $g^{-1} \in A$ and $d(g,A) = i_2(g)$ if $g^{-1}
  \notin A$. Since $\mathcal U(G) \cong \coprod_{G \cdot A \in G
  \backslash \mathcal U(G)} G \cdot A$, we have that
  \begin{eqnarray*}
  \lefteqn{  S_\bullet(\pi^G_e)(s+t)}\\
   &=& S_\bullet(\pi^G_e) S_+(\nabla)(s,t) \\
    &=& S[G/e \amalg G/e \xleftarrow d G/e \times \mathcal U(G)
  \xrightarrow{} \mathcal U(G) \to G/G](s,t) \\
  &=& \sum_{G \cdot A \in G
  \backslash \mathcal U(G)} S[G/e \amalg G/e \xleftarrow {d} G/e \times
  GA \rightarrow  GA \to G/G](s,t) \\
  &=& \!\!\!\!\!\!\!\! \sum_{G \cdot A \in G
  \backslash \mathcal U(G)} \!\!\!\!\!\!\!\! S[G/e \! \amalg \! G/e
  \xleftarrow {} G/e \! \times \!   
  A/U_A \! \amalg \! G/e \! \times \! (G - A)/U_A 
  \rightarrow  G/U_A \to 
  G/G](s,t) \\
  &=& \sum_{G \cdot A \in G
  \backslash \mathcal U(G)} S[G/e \xleftarrow {=} G/e
  \xrightarrow{\pi^{U_A}_e} G/U_A \xrightarrow {\pi^G_{U_A}}
  G/G](s^{i_A} t^{i_{G-A}}) \\
  &=& \sum_{G \cdot A \in G
  \backslash \mathcal U(G)} S_+({\pi^G_{U_A}}) S_\bullet(
  \pi^{U_A}_e)(s^{i_A} t^{i_{G-A}}),  
  \end{eqnarray*}
  where the maps without labels are natural projections. 
\end{proof}
\begin{lem}\label{sumlemma2}
  Let $R$ be a commutative ring. With the notation of lemma
  \ref{sumlemma}, we have for every subset $A$ and subgroup $U$ of $G$
  and for every $s,t \in R$, that
  \begin{displaymath}
    (s+t)^{(G:U)} = \sum_{G\cdot A \in G \backslash \mathcal
  U(G)} |U(G/U_A)^U| \cdot (s^{i_A} \cdot t^{i_{G-A}})^{(U_A:U)}.
  \end{displaymath}
\end{lem}
\begin{proof}
  We compute:
  \begin{eqnarray*}
    (s+t)^{(G:U)} &=& \sum_{A \subseteq G/U} s^{|A|} t^{|G/U| -
    |A|} =  \sum_{A \in \mathcal U(G), \,U \leq U_A} s^{|A/U|}
  t^{|(G-A)/U|} \\
  &=& \sum_{G\cdot A \in G \backslash \mathcal
  U(G)} |(G/U_A)^U| \cdot (s^{i_A} \cdot t^{i_{G-A}})^{(U_A:U)}.
  \end{eqnarray*}
\end{proof}
The following lemma is a variation on \cite[lemma 3.2.5]{DS}, and the
proof we give here is essentially identical to the one given in
\cite{DS}.  
\begin{lem}\label{hardsum}
  For some $k \in \N$ let $V_1,\dots,V_k \leq G$ be a sequence of
  subgroups of $G$. Then for every subgroup $U \leq G$ and every
  sequence $\varepsilon_1,\dots,\varepsilon_k \in \{\pm 1\}$ there
  exists a polynomial $\xi_U
  = \xi^G_{(U;V_1,\dots,V_k;\varepsilon_1,\dots,\varepsilon_k)} =
  \xi_U(x_1,\dots,x_k) \in \Z[x_1,\dots,x_k]$ such that for every
  $G$-Tambara functor $S$ and all $s_1,\dots,s_k \in S(G/e)$ one has
  \begin{displaymath}
    \sum_{i=1}^k \varepsilon_i S_+(\pi^G_{V_i}) S_\bullet(\pi^{V_i}_e)
  (s_i) = {\sum}'_{U \leq G} S_+(\pi^G_U) S_\bullet(\pi^U_e)(
  \xi_U(s_1,\dots,s_k)).
  \end{displaymath}
\end{lem}
\begin{proof}
  If $\varepsilon_1 = \varepsilon_2 = \dots = \varepsilon_k = 1$ and if
  $V_i$ is not conjugate to $V_j$ for $i \ne j$, then  
  \begin{displaymath}
    \sum_{i=1}^k \varepsilon_i S_+(\pi^G_{V_i}) S_\bullet(\pi^{V_i}_e)
  (s_i) = {\sum}'_{U \leq G} S_+(\pi^G_U) S_\bullet(S^U_e)(s_U), 
  \end{displaymath}
  with $s_U = s_i$ if $U$ is conjugate to $V_i$ and $s_U = 0$ if $U$ is
  not conjugate to any of the $V_1,\dots,V_k$. So in this case we are
  done: $  \xi_U(s_1,\dots,s_k) = s_U$. We use triple induction to prove
  the lemma. First with respect to $m_1 =
  m_1(U;V_1,\dots,V_k;\varepsilon_1,\dots,\varepsilon_k) := \max (
  |V_i| \, | \, \varepsilon_i = -1 \text{ or there exists some } j \ne i
  \text{ with } V_j \sim V_i)$, then with respect to 
  $m_3 := | \{ i \, | \, |V_i| = m_1 \text{
  and there exists some } j \ne i \text{ with } V_j \sim V_i\}|$
  and then
  with respect to 
  $m_2 := | \{ i
  \, | \, |V_i| = m_1 \text{ and } \varepsilon_i = -1\}|$.

  We have just verified that the lemma holds in the case $m_1 = 0$. In case
  $m_1 > 0$ we have either $m_2 > 0$ or $m_3 > 0$. In case $m_2 > 0$,
  say $|V_1| = m_1$ and $\varepsilon_1 = -1$, we may use
  \ref{sumlemma} with $G=V_1$, $s = -s_1$, $t = s_1$ to conclude that 
  \begin{displaymath}
  0 = \sum_{V_1 \cdot A \in V_1 \backslash \mathcal U(V_1)}
  S_+(\pi^{V_1}_{U_A}) S_\bullet(\pi^{U_A}_e) ((-1)^{i_A} s_1^{(V_1:U_A)}).
  \end{displaymath}
  Therefore, considering the two special summands $A = \emptyset$ and $A
  = V_1$ and putting $\mathcal U_0(V_1) := \{ A \in \mathcal U(V_1) \, |
  \, A \ne \emptyset \text{ and } A \ne V_1 \}$, one gets
  \begin{displaymath}
  -S_{\bullet}(\pi^{V_1}_e)(s_1) = S_{\bullet}(\pi^{V_1}_e)(-s_1) +
  \sum_{V_1 \cdot A \in V_1 \backslash \mathcal U_0(V_1)} 
  S_+(\pi^{V_1}_{U_A}) S_\bullet(\pi^{U_A}_e) ((-1)^{i_A} s_1^{(V_1:U_A)}).
  \end{displaymath}
  Hence, if $A_{k+1}, A_{k+2}, \dots, A_{k'} \in \mathcal U_0(V_1)$
  denote representatives of the $V_1$-orbits $V_1 A \subseteq \mathcal
  U_0(V_1)$ and we let $V_{k+1} := U_{A_{k+1}}, \dots, V_{k'} :=
  U_{A_{k'}}$ then $V_i \leq V_1$  and $V \ne V_1$ for $i \ge k+1$. If we put
  $\varepsilon_{k+1} = \dots = \varepsilon_{k'} = 
  1$ and $s_{k+1} := (-1)^{i_{A_{k+1}}} s_1^{(V_1:V_{k+1})}, \dots,
  s_{k'} = (-1)^{i_{A_{k'}}} s_1^{(V_1:V_{k'})}$, then the polynomial
  \begin{displaymath}
    \xi^G_{(U;V_1,\dots,V_{k};{-1},\varepsilon_2,\dots,\varepsilon_k)
  (s_1,\dots,s_k)} := 
    \xi^G_{(U;V_1,\dots,V_{k'};1,\varepsilon_2,\dots,\varepsilon_{k'})
  (-s_1,s_2,\dots,s_{k'})}
  \end{displaymath}
  makes the statement of the lemma hold.
  We can conclude that if the lemma holds for every $(n_1,n_2,n_3)$ with
  either $n_1 < m_1$ or ($n_1 = m_1$, $n_2 < m_2$ and $n_3 \leq m_3$),
  then it also holds for $(m_1,m_2,m_3)$.

  Similarly, if $m_2 = 0$, but $m_3 > 0$, say $V_1 = V_2$, then we may
  use lemma \ref{sumlemma} once more with $G = V_1$, $s = s_1$, and $t =
  s_2$ to conclude that
  \begin{eqnarray*}
  S_{\bullet}(\pi^{V_1}_e)(s_1 + s_2) &=& 
  \sum_{V_1 \cdot A \in V_1 \backslash \mathcal U(V_1)} 
  S_+(\pi^{V_1}_{U_A}) S_\bullet(\pi^{U_A}_e) (s_1^{i_A}
  s_2^{i_{V_1-A}}) \\
  &=& 
  S_{\bullet}(\pi^{V_1}_e)(s_1) + 
  S_{\bullet}(\pi^{V_1}_e)(s_2)  \\
  && \mbox{} +
  \sum_{V_1 \cdot A \in V_1 \backslash \mathcal U_0(V_1)} 
  S_+(\pi^{V_1}_{U_A}) S_\bullet(\pi^{U_A}_e) (s_1^{i_A}
  s_2^{i_{V_1-A}}),  
  \end{eqnarray*}
  so with $V_{k+1}, \dots,V_{k'}$ as above, but with $\varepsilon_{k+1} =
  \dots = \varepsilon_{k'} = -1$ and with $s_{k+1} := s_1^{i_{A_{k+1}}}
  s_2^{i_{V_1 - A_{k+1}}}, \dots, s_{k'} := s_1^{i_{A_{k'}}}
  s_2^{i_{V_1 - A_{k'}}}$, the polynomial
  \begin{displaymath}
    \xi^G_{(U;V_1,\dots,V_{k};\varepsilon_1,\dots,\varepsilon_k)
  (s_1,\dots,s_k)} := 
    \xi^G_{(U;V_2,\dots,V_{k'};1,\varepsilon_2,\dots,\varepsilon_{k'})
  (s_1+s_2,s_3,\dots,s_{k'})}
  \end{displaymath}
  makes the statement of the lemma hold. 
  We can conclude that if the lemma holds for every $(n_1,n_2,n_3)$ with
  either $n_1 < m_1$ or ($n_1=m_1$, $n_2 = m_2 = 0$ and $n_3 < m_3$), then
  it also holds for $(m_1,0,m_3)$. The statement of
  the lemma now follows by induction first on $m_1$, then on $m_3$ and
  finally on $m_2$.
\end{proof}
\begin{lem}\label{hardsum2}
  Let $R$ be a ring.
  For some $k \in \N$ let $V_1,\dots,V_k \leq G$ be a sequence of
  subgroups of $G$ and let $\varepsilon_1,\dots,\varepsilon_k \in \{\pm
  1\}$. For every subgroup $H$ of $G$ one has, with the notation of
  lemma \ref{hardsum}, that
  \begin{displaymath}
    \sum_{i=1}^k \varepsilon_i |(G/{V_i})^H| (s_i)^{(V_i:H)} =
  {\sum}'_{U \leq G} |(G/U)^H| (
  \xi_U(s_1,\dots,s_k))^{(U:H)}.
  \end{displaymath}
\end{lem}
\begin{proof}
  We use the same notation as in the above proof, and prove the lemma by
  induction on first $m_1$, then on $m_3$ and finally on $m_2$. Suppose
  that the lemma holds for every $(n_1,n_2,n_3)$ with $n_1 < m_1$ or
  $n_3 < m_3$ or ($n_1 = m_1$ and $n_3 = m_3$ and $n_2 < m_2$), and that
  we are given 
  $(V_1,\dots,V_k;\varepsilon_1,\dots,\varepsilon_k)$ with
  $m_i = m_i(V_1,\dots,V_k;\varepsilon_1,\dots,\varepsilon_k)$. In the
  above proof we constructed 
  $(W_1,\dots,W_{k'};\delta_1,\dots,\delta_{k'})$ and $t_1,\dots,t_{k'}
  \in R$ such that
  $n_i = n_i(W_1,\dots,W_{k'};\delta_1,\dots,\delta_{k'})$, $i =1,2,3$
  are on the form covered by the induction hypothesis, and such that:
  \begin{eqnarray*}
    \sum_{i=1}^k \varepsilon_i |(G/V_i)^H| s_i^{(V_i:H)} &=& 
    \sum_{i=1}^{k'} \delta_i |(G/W_i)^H| t_i^{(W_i:H)} \\
    &=& {\sum_{U \leq G}}' |(G/U)^H|
    (\xi^G_{(U;W_1,\dots,W_{k'};\delta_1,\dots,\delta_{k'})}
    (t_1,\dots,t_{k'}))^{(U:H)}  \\
    &=& {\sum_{U \leq G}}' |(G/U)^H|
    (\xi^G_{(U;V_1,\dots,V_{k};\varepsilon_1,\dots,\varepsilon_{k})}
    (s_1,\dots,s_{k}))^{(U:H)}.  
  \end{eqnarray*}
  The other cases of the induction are similar.
\end{proof}
\begin{lem}\label{prodlemma}
  For every $G$-Tambara functor $S$, for every two subgroups $V, W \leq
  G$ and for every $s,t \in S(G/e)$ one has the modified Mackey formula
  \begin{multline*}
    \lefteqn{S_+(\pi^G_V)S_\bullet(\pi^V_e)(s) \cdot 
    S_+(\pi^G_W)S_\bullet(\pi^W_e)(t)} \\ =
    \sum_{VgW \in V \backslash G /W} S_+(\pi^G_{V \cap gWg^{-1}})
  S_\bullet(\pi^{V \cap gWg^{-1}}_e) (s^{(V: V \cap gWg^{-1})} \cdot
  t^{(W:g^{-1}Vg \cap W)}).
  \end{multline*}
\end{lem}
\begin{proof}
  We use the diagram
  \begin{displaymath}
  \begin{CD}
   @. G/e \amalg G/e @<{\pr_1 \amalg \pr_2}<< G/e \times G/W \amalg G/V
  \times G/e \\ 
   @. @V{\pi^V_e \amalg \pi^W_e}VV @VV{\pi^V_e \times \id \amalg \id
  \times \pi^W_e}V \\
   G/G \amalg G/G @<{\pi^G_V \amalg \pi^G_W}<< G/V \amalg G/W @<{\pr_1
  \amalg \pr_2}<< G/V 
  \times G/W \amalg G/V 
  \times G/W \\
   @V{\nabla}VV @. @V{\nabla}VV \\
   G/G @<<< G/V \times G/W @= G/V \times G/W,
  \end{CD}
  \end{displaymath}
  where $\nabla$ is the fold map, the upper square is a
  pull-back and the lower rectangle is an exponential diagram.
  Concatenating with the diagram 
  \begin{displaymath}
  \begin{CD}
   G/e \times G/W \amalg G/V \times G/e @<{\alpha_1 \amalg
  \alpha_2}<{\cong}< \coprod_{gW \in G/W} G/e 
  \amalg \coprod_{V\! g \in V \backslash G} G/e \\
   @V{\pi^V_e \times \id \amalg \id \times \pi^W_e}VV @| \\
   G/V \times G/W \amalg G/V \times G/W 
   @. \coprod_{gW \in G/W} G/e 
   \amalg \coprod_{V\! g \in V \backslash G} G/e \\
   @V{\nabla}VV @V{(\delta_1,\delta_2)}VV \\
   G/V \times G/W @<{\gamma}<{\cong}< \coprod_{VgW \in V \backslash G /
  W} G/(V \cap gWg^{-1}) 
  \end{CD}
  \end{displaymath}
  with maps defined by 
  \begin{eqnarray*}
   \alpha_1(gW, \sigma) &=& (\sigma,\sigma g W), \\
   \alpha_2(Vg,\sigma) &=& (\sigma g^{-1}V, \sigma), \\
   \delta_1(gW,\sigma) &=& (VgW,\sigma(V \cap gWg^{-1})), \\ 
   \delta_2(Vg,\sigma) &=& (VgW, \sigma g^{-1}(V \cap gWg^{-1})) \qquad
  \text{and} \\ 
  \gamma(VgW,\sigma(V \cap gWg^{-1})) &=& (\sigma V,\sigma g W),
  \end{eqnarray*}
   we get that
  \begin{eqnarray*}
    \lefteqn{
    S_+(\pi^G_V)S_\bullet(\pi^V_e)(s) \cdot 
    S_+(\pi^G_W)S_\bullet(\pi^W_e)(t)
    } \\
    &=& 
    S[G/e \amalg G/e \leftarrow \!\!\!\! 
  \coprod_{gW \in G/W} 
  \!\!\!\!  G/e
  \amalg \!\!\!\! 
  \coprod_{Vg \in V \backslash G}
  \!\!\!\! G/e \to \!\!\!\!
  \!\!\!\! \coprod_{VgW \in 
  V \backslash G 
  /W} \!\!\!\! \!\!\!\! G/(V \cap gWg^{-1}) \to G/G](s,t) \\
  &=&
  \sum_{VgW \in V \backslash G / W} S_+(\pi^G_{V \cap gWg^{-1}})
  S_\bullet(\pi^{V \cap gWg^{-1}}_e) (s^{(V:{V \cap gWg^{-1}})} \cdot
  t^{(W:g^{-1}Vg \cap W)}).
  \end{eqnarray*}
\end{proof}
\begin{lem}[{\cite[Lemma 3.2.13]{DS}}]
\label{prodlemma2}
  For every ring $R$, for every three subgroups $H, V, W \leq
  G$ and for every $s,t \in R$ one has the modified Mackey formula
  \begin{multline*}
    \lefteqn{|(G/V)^H|s^{(V:H)} |(G/W)^H| t^{(W:H)}} \\ =
    \sum_{VgW \in V \backslash G /W} |(G/{V \cap gWg^{-1}})^H|
    (s^{(V: V \cap gWg^{-1})} \cdot
    t^{(W:g^{-1}Vg \cap W)})^{({V \cap gWg^{-1}}:H)}.
  \end{multline*}
\end{lem}
\begin{proof}[Proof of theorem \ref{maintam}]
  Let $G= V_1,V_2,\dots,V_k = U$ be a system of representatives of
  subgroups of $G$ containing a conjugate of $U$. We define 
  \begin{displaymath}
    s^G_U(a_{V_1},b_{V_1}, \dots,a_{V_k},b_{V_k}) :=
  \xi^G_{(U;V_1,V_1,\dots,V_k,V_k;1,\dots,1)} (a_{V_1},b_{V_1},
  \dots,a_{V_k},b_{V_k}) 
  \end{displaymath}
  and
  \begin{displaymath}
    m^G_U(a_{V_1}, \dots,a_{V_k}) :=
    \xi^G_{(U;V_1,\dots,V_k;-1,\dots,-1)}
    (a_{V_1}, \dots,a_{V_k}).
  \end{displaymath}
  By the lemmas \ref{hardsum} and \ref{hardsum2} these are integral
  polynomials with the desired 
  properties. For example, for $U \leq
  G$ we have that 
  \begin{eqnarray*}
    \phi_U(a)+\phi_U(b) &=& \sum_{i=1}^k |(G/V_i)^U| (a_{V_i}^{(V_i:U)} +
  b_{V_i}^{(V_i:U)}) \\
    &=&
  \phi_U(\xi_{(U;V_1,V_1,\dots,V_k,V_k;1,\dots,1)}(a_{V_1},b_{V_1},
  \dots,a_{V_k}, b_{V_k})) \\ &=& \phi_U(s_U(a_{V_1},b_{V_1},
  \dots,a_{V_k}, b_{V_k})).
  \end{eqnarray*}
  To construct 
  $p_U = p^G_U$ we first choose a system
  $x_1,x_2 \dots, x_h$ of representatives of the $G$-orbits in 
  \begin{displaymath}
    X := \coprod_{i,j = 1}^k G/V_i \times G/V_j.
  \end{displaymath}
  Next we put $W_r := G_{x_r}$ and
  \begin{displaymath}
    p_r = p_r(a_{V_1},b_{V_1},\dots,a_{V_k},b_{V_k}) := a_i^{(V_i:W_r)} \cdot
  b_j^{(V_j:W_r)} 
  \end{displaymath}
  in case $x_r = (g_r V_i,g'_rV_j) \in G/V_i \times G/V_j \subseteq
  X$. Using these conventions, we define
  \begin{displaymath}
    p_U^G(a_{V_1},b_{V_1},\dots,a_{V_k},b_{V_k}) := 
    \xi^G_{(U;W_1,\dots,W_h;1,\dots,1)}(p_1,\dots,p_r).
  \end{displaymath}
  Using the lemmas \ref{prodlemma} and \ref{prodlemma2} we see that $p_U$ has the desired
  properties. 
\end{proof}
\section{Free Tambara functors}
\label{freeta}
In this section we give a proof of theorem \ref{mainwitt}.
On the way we shall give a combinatorial description of the
Witt vectors of a polynomial $G$-ring, that is, a $G$-ring of the form
$U^G(X,G/e)$ 
for a finite $G$-set $X$.
Recall from \ref{Gringistam} that the functor $R \mapsto R(G/e)$ from
the category $[U^{fG},\Set]_0$ of $fG$-Tambara functors 
to the category of $G$-rings is an equivalence of
categories, and that there 
are morphisms $(\Ens^{fG})^{\op} \subseteq (\Ens^G)^{\op} \xrightarrow R
U^G$ of $\uO(G)$-coloured theories. We let $F : \Set^G \simeq
[(\Ens^G)^{\op},\Set]_0 \to [U^G,\Set]_0$ denote the left adjoint of the
forgetful functor induced by the above composition of morphisms of
coloured theories.
\begin{defn}
  Given finite $G$-sets $X$ and $Y$ we let $\widetilde U^G_+(X,Y) \subseteq
  \overline U^G(X,Y)$ denote those elements of the form $X \leftarrow A
  \to B \to 
  Y$, where $G$ acts freely on $A$, and we let $\widetilde U^G(X,Y)
  \subseteq U^G(X,Y)$ denote the abelian subgroup generated by
  $\widetilde U^G_+(X,Y)$. 
  The composition 
  \begin{displaymath}
    U^G(Y,Z) \!\! \times \!\! \widetilde U^G(X,Y) \cup \widetilde U^G(Y,Z) \!\!
  \times \!\! U^G(X,Y) \subseteq
  U^G(Y,Z) \!\! \times \!\! U^G(X,Y) 
  \xrightarrow{\circ} U^G(X,Z)
  \end{displaymath}
  factors through the inclusion $\widetilde U^G(X,Z) \subseteq
  U^G(X,Z)$. 
  We obtain a functor 
  $\widetilde U^G : \Ens^G \to [U^G,\Set]_0$ with
  $\widetilde U^G(f : Y \to X) = \widetilde U^G(R_f,-)$. 
\end{defn}
\begin{lem}\label{tamGset}
  Given a $G$-Tambara functor $S$ and a finite free $G$-set $A$, there is an
  isomorphism $\Set^G(A,S(G/e)) \xrightarrow \cong S^*(A)$, which is
  natural in $A$.
\end{lem}
\begin{proof}
  Choosing an isomorphism $\phi :A \xrightarrow \cong G/e \times A_0$
  we obtain an isomorphism
  \begin{eqnarray*}
    \Set^G(A,S(G/e)) &\xrightarrow {{\phi^{-1}}^*}& \Set^G(G/e \times A_0,
    S(G/e)) \\ 
    &\cong& \Set(A_0,S(G/e)) \cong S(G/e \times A_0)
    \xrightarrow {S^*(\phi)} S(A).
  \end{eqnarray*}
  This isomorphism is independent of the choise of $\phi$.
\end{proof}
\begin{lem}
  Given a finite $G$-set $X$ we have that $FX \cong \widetilde U^G(X,-)$.
\end{lem}
\begin{proof}
  For every $G$-Tambara functor $S$ we shall construct a bijection
  $$\Set^G(X,S(G/e)) \cong [U^G,\Set]_0(\widetilde U^G(X,-),S).$$ 
  Given $f
  : X \to S(G/e) \in \Set^G(X,S(G/e))$ we let $\phi(f) \in
  [U^G,\Set]_0(\widetilde U^G(X,-),S)$ take $x = [X \xleftarrow d
  A \xrightarrow b B \xrightarrow c Y] \in \widetilde U^G(X,Y)$ to
  $\phi(f)(x) \in S(Y)$ constructed as follows:  
  by \ref{tamGset} we obtain an 
  element $a \in S(A)$, and we let $\phi(f)(x) =
  S_+(c)S_\bullet(b)(a)$. Conversely, given $g \in
  [U^G,\Set]_0(\widetilde U^G(X,x),S)$, we construct $\psi(g) \in 
  \Set^G(X,S(G/e))$ by letting $\psi(g)(x) =  [X \leftarrow
  G/e \to G/e \to G/e]$, where the map pointing left takes $e \in G$ to
  $x \in X$
  and where the maps pointing right are identity maps. We leave it to
  the reader to check that $\phi$ and $\psi$ are inverse bijections. 
\end{proof}
\begin{cor}
  For every finite $G$-set $X$ the functor $\widetilde U^G(X,-) : U^G \to \Set$
  is isomorphic to $L_GU^G(X,G/e)$.
\end{cor}
\begin{thm}\label{subwitt}
  Let $X$ be a finite $G$-set and let $R = \widetilde U^G(X,-)$. Then
  $R(G/e) = 
  U^G(X,G/e)$, and the Teichm\"uller
  homomorphism $$\tau :
  \W_G(L_GR(G/e))/{\IW}_G(L_GR(G/e)) \to L_GR(G/G) = \widetilde
  U^G(X,G/G)$$
  is an isomorphism.
\end{thm}
\begin{proof}[Proof of theorem \ref{mainwitt0}]
  Let $X$ be a set with trivial action of $G$. The polynomial ring
  $\Z[X]$ has the universal property that determines the left adjoint of
  the forgetful functor from rings with an action of $G$ to $\Set^G$ 
  up to isomorphism. Hence $\Z[X] \cong \widetilde U^G(X,G/e) =
  U^G(X,G/e)$. By theorem \ref{subwitt} $\widetilde U^G(X,G/G) =
  U^G(X,G/G)$ is isomorphic to $\W_G(\Z[X])$.
\end{proof}
\begin{proof}[Proof of theorem \ref{mainwitt}]
  Let $R_0 = R(G/e)$.
  Given $\alpha = [W \leftarrow C \to D \to X] \in U^{fG}(W,X)$ we
  have an $fG$-Tambara map $\alpha^* : U^G(X,-) \to U^G(W,-)$ and we
  have the map $R(\alpha) : R(W) \to R(X)$. Hence we obtain maps
  \begin{displaymath}
      U^G(X,G/G) \times R(X) \leftarrow U^G(X,G/G) \times R(W) \to
      U^G(W,G/G) \times R(W).
  \end{displaymath}
  Recall from lemma \ref{leftKan} that $L_GR(G/G)$ can be constructed by
  the coequalizer 
  of the diagram
  \begin{displaymath}
    \coprod_{X,Y \in \ob U^{fG}} U^G(X,G/G) \times R(Y) \stackrel
    {\displaystyle \rightarrow} \rightarrow \coprod_{X \in \ob U^{fG}}
  U^G(X,G/G) 
    \times R(X),
  \end{displaymath}
  induced by the above maps.
  We shall construct a map $$\rho : L_GR(G/G) \to \W_G(R_0)/{\IW}_G(R_0)$$
  by specifying explicit maps $\rho_X : U^G(X,G/G) \times R(X) \to
  \W_G(R_0)/{\IW}_G(R_0)$. Given $\underline r \in R(X)$, we have an
  $fG$-Tambara morphism $\ev_{\underline r} : U^{fG}(X,-) \to R$. In
  particular we have a $G$-homomorphism $\ev_{\underline r}(G/e) : U^G(X,G/e) \to R(G/e) =
  R_0$. By theorem \ref{subwitt} we get an induced
  ring-homomorphism
  \begin{displaymath}
    U^G(X,G/G) \cong \W_G(U^G(X,G/e))/{\IW}_G(U^G(X,G/e)) \to \W_G(R_0)/{\IW}_G(R_0),
  \end{displaymath}
  and hence by adjunction we obtain a map
  \begin{displaymath}
    \rho_X : U^G(X,G/G) \times R(X) \to \W_G(R_0)/{\IW}_G(R_0).
  \end{displaymath}
  We need to check that these $\rho_X$ induce a map on the coequalizer
  $L_GR(G/G)$ of the above coequalizer diagram, that is, for $\alpha$
  as above we need to show that the diagram 
  \begin{displaymath}
    \begin{CD}
      U^G(X,G/G) \times R(W) @>{\alpha^* \times \id}>>
      U^G(W,G/G) \times R(W) \\
      @V{\id \times R(\alpha)}VV @V{\rho_W}VV \\
      U^G(X,G/G) \times R(X) @>{\rho_X}>> \W_G(R_0)/{\IW}_G(R_0)       
    \end{CD}
  \end{displaymath}
  commutes. In order to do this we note that the diagram
  \begin{displaymath}
    \begin{CD}
      \W_G(U^G(X,G/e)) @>{t}>> U^G(X,G/G) \\
      @V{\alpha^*}VV @V{\alpha^*}VV \\
      \W_G(U^G(W,G/e)) @>{t}>> U^G(W,G/G) 
    \end{CD}
  \end{displaymath}
  commutes, and therefore it will suffice to show that the diagram
  \begin{displaymath}
    \begin{CD}
      \W_G(U^G(X,G/e)) \times R(W) @>{\alpha^* \times \id}>> 
      \W_G(U^G(W,G/e)) \times R(W) \\
       @V{\id \times R(\alpha)}VV @VVV \\
      \W_G(U^G(X,G/e)) \times R(X) @>>> \W_G(R_0)
    \end{CD}
  \end{displaymath}
  commutes, where the arrows without labels are constructed using 
  the homomorphisms $\W_G(\ev_{\underline r}(G/e))$ for $\underline r$
  an element 
  of either $R(X)$ or $R(W)$. Using diagonal inclusions
  of the form 
  \begin{displaymath}
    \W_G(T) \times Z \rightarrow \W_G(T) \times {\prod_{U \leq G}}' Z \approx
    {\prod_{U \leq G}}'(T \times Z) 
  \end{displaymath}
  we see that it suffices to show note that the diagram
  \begin{displaymath}
    \begin{CD}
      {\prod}'_{U \leq G} (U^G(X,G/e) \times R(W))
      @>{{\prod}'_{U \leq G} (\alpha^* \times \id)}>>
      {\prod}'_{U \leq G} (U^G(W,G/e) \times R(W)) \\
      @V{{\prod}'_{U \leq G} (\id \times R(\alpha))}VV @VVV \\      
      {\prod}'_{U \leq G} (U^G(X,G/e) \times R(X))
      @>>> 
      {\prod}'_{U \leq G} R(G/e) \approx \W_G(R_0)
    \end{CD}
  \end{displaymath}
  commutes. This finishes the construction of $\rho : L_GR(G/G) \to
  \W_G(R_0)/{\IW}_G(R_0)$. We leave it to the reader to check that $\rho$
  and $\tau$ are inverse bijections. For this it might be helpful to
  note that
  \begin{eqnarray*}
   && {\sum_{U \leq G}}' ([G/e
    \xleftarrow = G/e \xrightarrow {\pi^U_e} G/U \rightarrow G/G]
    \circ [Y \leftarrow A_U \to B_U \to G/e]) \\
    &=& [\!{\coprod_{U \leq G}}' \! G/e \xleftarrow = \!\! {\coprod_{U
  \leq G}}' \!
    G/e \! \to \! \! {\coprod_{U \leq G}}' \! G/U \to G/G] \! \circ \!
  [Y \! \leftarrow \!
    \! {\coprod_{U \leq G}}' \! A_U \! \to \! {\coprod_{U \leq G}}' \!
  B_U \! \to
    \! {\coprod_{U \leq G}}'  \! G/e]. 
  \end{eqnarray*}
\end{proof}
In order to begin the proof of theorem \ref{subwitt} we need to introduce
filtrations of both sides.
\begin{defn}
  Let $R$ be a ring and let $U \leq G$ be a subgroup of $G$. We let
  $I_U \subseteq \W_G(R)$ denote the ideal generated by those $a =
  (a_K)'_{K \leq G} \in \W_G(R)$ satisfying that $a_K \ne 0$
  implies that $K \seq U$. We let $\widetilde I_U \subseteq I_U$
  denote the subideal $\widetilde I_U = \sum_{V \subsetneq U} I_V
  \subseteq I_U$. 
\end{defn}
\begin{defn}
  Given a $G$-sets $X$, we let 
  $J^+_U = J_U(X,G/G) \subseteq \widetilde U^G(X,G/G)$ denote
  the subset of elements of the form $[X \leftarrow A \to B \to
  G/G] \in U^G_+(X,G/G) \subseteq U^G(X,G/G)$, where $B^K = \emptyset$ when $U
  \seq K$ and $U \ne K$.
  We let $J_U \subseteq U^G(X,G/G)$ denote the ideal
  generated by $J_U^+$ and we let $\widetilde J_U \subseteq J_U$
  denote the subideal $\widetilde J_U = \sum_{V \subsetneq U} J_V
  \subseteq J_U$. 
\end{defn}
\begin{lem}\label{omnibus}
  In the situation of the above definition the following holds:
  \begin{enumerate}
    \item[(i)] Any element in $J_U$ is of the form $x-y$ for $x,y \in  J_U^+$.
    \item[(ii)] Every element $x$ in the image of the map $J^+_U \to J_U/
      \widetilde J_U$ is of the form
      \begin{displaymath}
        x = [X \xleftarrow d G/e \times A \xrightarrow {\pi^U_e \times f} G/U
        \times B \xrightarrow q G/G] + \widetilde J_U, 
      \end{displaymath}
      with $f: A \to B$ a map of (nonequivariant) sets and $d$ a $G$-map, 
      where $q$ is the composition $G/U \times B
      \xrightarrow{\pr} G/U \xrightarrow {\pi^G_U} G/G$.
    \item[(iii)] If 
      \begin{displaymath}
        x = [X \xleftarrow d G/e \times A \xrightarrow {\pi^U_e \times f} G/U
        \times B \xrightarrow q G/G] + \widetilde J_U, 
      \end{displaymath}
      and
      \begin{displaymath}
        x' = [X \xleftarrow {d'} G/e \times A' \xrightarrow {\pi^U_e \times f'} G/U
        \times B' \xrightarrow {q'} G/G] + \widetilde J_U, 
      \end{displaymath}
      then $x = x'$ if and only if there exist
      bijections $\alpha : A \to A'$ and $\beta : B \to B'$ and for
      every $a \in A$ there exists $g_a \in N_G(U)$ such that
      firstly $d'(e,\alpha a) = d(g_a,a)$, and secondly, if $a_1,a_2 \in A$
      satisfy that $f(a_1) = f(a_2)$, then $g_{a_1} U = g_{a_2}
      U$. 
  \end{enumerate}
\end{lem}
\begin{proof}
  A straight forward verification yields that multiplication in
  $\widetilde U^G_+(X,G/G)$ induces a map $J^+_U \times \widetilde U^G_+(X,G/G)
  \to J^+_U$ and that $J_U^+$ is closed under sum. It follows that
  $J_U$ is the abelian subgroup of 
  $U^G(X,G/G)$ generated by $J^+_U$. The statement (i) is a direct
  consequence of this.
  For (ii) we note that for every element 
      \begin{displaymath}
        r = [X \xleftarrow d D \xrightarrow {c} E \xrightarrow t G/G]  
      \end{displaymath}
  in $\widetilde U^G_+(X,G/G)$, we have a decomposition $E \cong
  {\coprod}'_{K\leq G} E_K$, where $E_K \cong G/K \times B_K$ for
  some $B_K$. This decomposition induces an isomorphism
  \begin{displaymath}
    \widetilde U^G(X,G/G) \cong {\bigoplus_{K\leq G}}' J_K/
    \widetilde J_K
  \end{displaymath}
  of abelian groups. 
  Given an element $x$ of the form
      \begin{displaymath}
        x = [X \xleftarrow d D \xrightarrow {c} G/U
        \times B \xrightarrow q G/G] + \widetilde J_U, 
      \end{displaymath}
  we can choose a $G$-bijection of the form $c^{-1}(G/U \times \{b
  \}) \cong G/e \times A_b$ for every $b \in B$. It follows that $x$
  is represented by an element of the form 
      \begin{displaymath}
        x = [X \xleftarrow d G/e \times A \xrightarrow {p \times f} G/U
        \times B \xrightarrow q G/G] + \widetilde J_U.
      \end{displaymath}
  We leave the straight forward verification of part (iii) to the reader.
\end{proof}
\begin{lem}\label{prodrel}
  Let $U$ be a subgroup of $G$ and let
  $a = (a_U)'_{U \leq G} \in I_U$ have
  $$a_U = [X \xleftarrow  d A \times G/e \overset{f
  \times 1} \to B \times G/e \xrightarrow {\pr} G/e].$$ 
  Then 
  $$\tau(a) \equiv [X
  \xleftarrow d A \times G/e \overset{f 
  \times \pi^U_e} \to B \times G/U \rightarrow G/G] \quad \text{mod }
  \widetilde J_U.$$
\end{lem}
\begin{proof}
  The lemma follows from the diagram
  \begin{displaymath}
  \begin{CD}
    X @<d<< A \times G/e @<<< W @<<< G/e \times A \\
    @. @V{f \times \id}VV @VVV @VVV \\
    G/e @<{\pr}<< G/e \times B @<<< Z @<<< G/e \times B \\
    @V{\pi^U_e}VV @. @VVV @VVV \\
    G/U @<<< \coprod_{gU \in G/U} \map(gU,B) @= Y @<<< G/U \times B, 
  \end{CD}
  \end{displaymath}
  where the lower rectangle is an exponential diagram and the squares
  are pull-backs. We use that the map $G/U \times B \to Y$ which takes
  $(gU,b)$ to the 
  constant map $gU \to B$ with value $b$ is an isomorphism on
  $G/U$-parts and that $Y^H = \emptyset$ for $U \seq H$  and $U \ne H$.
\end{proof}
\begin{cor}
  Let $X$ be a $G$-set and let $R = \widetilde U^G(X,-)$.
  The map
  $$\tau : \W_G(U^G(X,G/e)) \to L_GR(G/G) = \widetilde U^G(X,G/G)$$
  satisfies that 
  $\tau(I_U) \subseteq 
  J_U$ and that $\tau(\widetilde I_U) \subseteq
  \widetilde J_U$. 
\end{cor}
\begin{prop}\label{gradediso}
  Let $X$ be a $G$-set, let $R = \widetilde U^G(X,-)$ and let $R_0 =
  R(G/e)$. For every $U 
  \leq G$ the map $\tau :
  \W_G(R_0)/{\IW}_G(R_0) \to L_GR(G/G) = \widetilde U^G(X,G/G)$ induces an
  isomorphism 
  $\tau_U : ({\IW}_G(R_0) + I_U)/({\IW}_G(R_0) + \widetilde I_U) \to J_U /
  \widetilde J_U$. 
\end{prop}
\begin{proof}
  Let $x \in I_U$ have $x_U = [X \xleftarrow  d A \times G/e \overset{f
  \times 1} \to B \times G/e \xrightarrow {\pr} G/e]$. Then by lemma
  \ref{prodrel} $\tau(x) \equiv [X
  \xleftarrow d A \times G/e \overset{f 
  \times p} \to B \times G/U \xrightarrow{q} G/G]$ mod $\widetilde
  J_U$, with the notation introduced there, and it follows from lemma
  \ref{omnibus} that $\tau_{U}$ is onto. 
  On the other hand, to prove injectivity, we pick $x_1,x_2\in I_U$ with
  $\tau(x_1) \equiv \tau(x_2) \mod \widetilde J_U$. Suppose that $x_{i,U}$
  has the form
  \begin{displaymath}
    x_{i,U} = [Z \xleftarrow {d_i} A_i \times G/e \overset{f_i
  \times 1} \to B_i \times G/e \xrightarrow = G/e]
  \end{displaymath}
  for $i = 1,2$.
  Let 
  \begin{displaymath}
    y_i = [Z \xleftarrow {d_i} A_i \times G/e \overset{f_i
  \times \pi^U_e} \to B \times G/U \xrightarrow q G/G]  
  \end{displaymath}
  for $i =1,2$.
  Then by lemma \ref{prodrel} $y_i \equiv \tau(x_i) \mod \widetilde J_U$
  for $i=1,2$.
  It follows from lemma \ref{omnibus} that
  there exist
      bijections $\alpha : A_1 \to A_2$ and $\beta : B_1 \to B_2$ and for
      every $a \in A_1$ there exists $g_a \in N_G(U)$ such that
      firstly $d_1(\alpha a) = g_a d(a)$ and secondly, if $a_1,a_2 \in A$
      satisfy that $f(a_1) = f(a_2)$, then $g_{a_1} U = g_{a_2}
      U$. 
  Given $a \in A_i$, let
  $z_{i,a} \in \widetilde U^G(X,G/e)$ denote the element $[X \xleftarrow
  {d_{i,a}} G/e 
  \xrightarrow p G/e \xrightarrow = G/e]$, where $d_{i,a}(e) =
  d_i(a,e)$. Then $z_{2,\alpha(a)} = g_a z_{1,a}$ and $x_{i,U} = \sum_{b \in
  B_i} ( \prod_{a \in f_i^{-1}(b)} 
  z_{i,a})$ for $i=1,2$, where an empty product is $1$ and an empty sum
  is $0$. We can conclude that $x_{1,U} - x_{2,U} \in {\IW}_G(R)$, and hence
  $x_1-x_2 \in {\IW}_G(R) + \widetilde I_U$.
  In the general case $\tau(x_1-x'_1) \equiv \tau(x_2-x'_2) \mod \widetilde
  J_U$ we 
  easily obtain that $x_1-x'_1 \equiv x_2-x'_2 \mod {\IW}_G(R) + \widetilde
  I_U$ by 
  collecting the positive terms.
\end{proof}
\begin{proof}[Proof of theorem \ref{subwitt}]
We first note that 
$\widetilde I_V =
\sum_{U \subsetneq V} I_U \cong
\colim{U \subsetneq V} I_U  
\subseteq I_V$ and that 
$\widetilde J_V = 
\sum_{U \subsetneq V} J_U \cong
\colim{U \subsetneq V} J_U 
\subseteq J_V$. The result now follows by induction on the
cardinality of $V$ using the above proposition and the five lemma on the
following map of
short exact 
sequences:
\begin{displaymath}
\begin{CD}
  \widetilde {\IW}_G(R(G/e)) + I_V @>>> {\IW}_G(R(G/e)) + I_V @>>>
  \frac{({\IW}_G(R(G/e)) + I_V)}{({\IW}_G(R(G/e)) + \widetilde I_V)} \\
  @VVV @VVV @VVV \\
  \widetilde J_V @>>> J_V @>>> J_V/\widetilde J_V.
\end{CD}
\end{displaymath}
\end{proof}

\section{The Witt Tambara-functor}
\label{witafun}
In this section we finally prove theorem \ref{mainwitt1}. Given a
subgroup $H \leq G$ and an $H$-set $X$, we can construct a $G \times
H$-set $G/e \times X$, where $G$ acts by multiplication on the left on
$G/e$, and where $h \cdot (g,x) := (gh^{-1},hx)$. We let $\ind^G_H X$
denote the $G$-set $H \backslash (G/e \times X)$.
\begin{lem}
  Let $H$ be a subgroup of $G$.
  The functor $\ind^G_H : \Ens^H \to \Ens^G $ induces functors
  $\ind^{G}_{H} : U^H \to U^G$, and  
  $\ind^{fG}_{fH} : U^{fH} \to U^{fG}$.
\end{lem}
\begin{proof}
  Since the functor $\ind^G_H : \Ens^H \to \Ens^G$ preserves
  pull-back diagrams and exponential diagrams it induces a functor 
  $\ind^G_H : U^H \to U^G$ that takes $X \leftarrow A \to B \to Y$
  to $\ind^G_H X \leftarrow \ind^G_H A \to \ind^G_H B \to \ind^G_H Y$. 
\end{proof}
Given a $G$-Tambara functor $S$ we can construct an $H$-Tambara functor
$\res^G_H S = S \circ \ind^G_H$. Similarly,  given an 
$fG$-Tambara functor $R$ we can construct an $fH$-Tambara functor
$\res^{fG}_{fH} S = S \circ \ind^{fG}_{fH}$. 
\begin{thm}\label{mainwitt2}
  Given an $fG$-Tambara functor $R$, the Teichm\"uller
  homomorphism 
  \begin{displaymath}  
    \tau : \W_H(\res^{fG}_{fH}R(H/e))/{\IW}_H(\res^{fG}_{fH}R(H/e)) 
    \rightarrow 
    \res^G_H L_G(H/H)
  \end{displaymath}
  is an 
  isomorphism. 
\end{thm}
\begin{proof}[Proof of theorem \ref{mainwitt1}]
If we consider $R_0 = R(G/e)$ as an $H$-ring, then
$$\W_H(\res^{fG}_{fH}(H/e))/{\IW}_H(\res^{fG}_{fH}(H/e)) \cong
\W_H(R_0)/{\IW}_H(R_0),$$  
and
$\res^G_H L_G(H/H) = L_GR(G/H)$. Combining these observations with
\ref{mainwitt2} we obtain the statement of theorem
\ref{mainwitt1}.
\end{proof}
\begin{lem}
  Let $H$ be a subgroup of $G$. The forgetful functor $i^*: \Ens^G  \to
  \Ens^H$
  which takes a $G$-set $Y$ to the same set considered as an $H$-set
  induces functors $i^* : U^G \to U^H$ and $i^*_f : U^{fG} \to U^{fH}$.
\end{lem}
\begin{lem}
  Let $H$ be a subgroup of $G$.
  The functor $i^* : U^G \to U^H$ is left adjoint to the functor
  $\ind^G_H : U^H \to U^G$ and 
  the functor $i^*_f : U^{fG} \to U^{fH}$ is left adjoint to the functor
  $\ind^{fG}_{fH} : U^{fH} \to U^{fG}$.
\end{lem}
\begin{proof}
  We prove only the first part of the lemma.
  Given $X \leftarrow A \to B \to G \times_H Y$ in $U^G(X,G \times_H Y)$
  we construct an element in $U^H(i^*X,Y)$ by the following diagram where
  the two squares furthest to the right are pull-back squares:
  \begin{displaymath}
  \begin{CD}
  i^*X @<<< A_H @>>> B_H @>>> H \times_H Y \\
  @VVV @VVV @VVV @VVV \\
  X @<<< A @>>> B @>>> G \times_H Y.
  \end{CD}
  \end{displaymath}
  Conversely, given $i^*X \leftarrow E \to F \to Y$ in $U^H(i^*X,Y)$ we
  construct the element $X \leftarrow G \times_H E \to G \times_H F \to
  G\times_H Y$ in $U^G(X,G \times_H Y)$. Here the arrow pointing to the
  left is the composite $G \times_H E \to G \times_H i^* X \to X$. We leave
  it to the reader to check that the maps are inverse bijections in an
  adjunction.
\end{proof}
We have the following commutative diagram of categories:
\begin{displaymath}
\begin{CD}
  U^{fG} @>{i^*_f}>> U^{fH} \\
  @VVV @VVV \\
  U^G @>{i^*}>> U^H,
\end{CD}
\end{displaymath}
where the vertical functors are the natural inclusions. 
Since $\ind^G_H$ is right adjoint to
$i^*$, $\res^G_H$ is left adjoint to
$[i^*,\Set]_0$ (see for example
\cite[Proposition 16.6.3]{Schubert}.) Similarly $\res^{fG}_{fH}$ is left 
adjoint to $[i_f^*,\Set]_0$.
From the commutative diagram of 
functor categories
\begin{displaymath}
\begin{CD}
  [U^{fG}, \Set]_0 @<{[i^*_f,\Set]_0}<< [U^{fH}, \Set]_0 \\
  @AAA @AAA \\
  [U^G, \Set]_0 @<{[i^*,\Set]_0}<< [U^H, \Set]_0,
\end{CD}
\end{displaymath}
where the vertical maps are the forgetful functors induced by the
inclusions $U^{fG} \subseteq U^G$ and $U^{fH} \subseteq U^H$ we can
conclude that 
there is a natural
isomorphism $\res^G_H L_G \cong L_H \res^{fG}_{fH}$.
\begin{proof}[Proof of theorem \ref{mainwitt2}]
  By theorem \ref{mainwitt} we have an isomorphism
  \begin{displaymath}
    \W_H(\res^{fG}_{fH} R(H/e))/
    {\IW}_H(\res^{fG}_{fH} R(H/e)) \xrightarrow \tau L_H \res^{fG}_{fH} R
  (H/H) \cong \res^G_H L_G R(H/H).
  \end{displaymath}
\end{proof}

\end{document}